




\input amstex
\documentstyle{amsppt}
\magnification\magstep1
\tolerance=2000
\NoBlackBoxes
\loadeusm
\loadeufm

\def\A{\eusm{A}}

\def\B{\eusm{B}}
\def\C{\eusm{C}}
\def\E{\Cal E}
\def\F{\eusm{F}(\eusm{H})}
\def\H{\eusm{H}}
\def\I{\eusm{I}}
\def\K{\eusm{K}}
\def\l{\lambda}
\def\la{\langle}
\def\N{\Bbb N}
\def\ot{\otimes}
\def\P{\Phi}

\def\ra{\rangle}
\def\rng{\operatorname{rng}}
\def\S{\Psi}
\def\tr{\operatorname{tr}}
\def\X{\eusm{X}}
\def\vp{\varphi}

\refstyle{A}
\widestnumber\key{ABCD}
\topmatter
\title
REFLEXIVITY OF THE AUTOMORPHISM AND ISOMETRY GROUPS OF
SOME STANDARD OPERATOR ALGEBRAS
\endtitle
\rightheadtext{Reflexivity of automorphism and isometry groups}
\author
LAJOS MOLN\' AR
\endauthor
\subjclass
Primary 47B49, 47D25, 46L40. Secondary 47B35
\endsubjclass
\keywords
Topological reflexivity, algebraic reflexivity, automorphism group,
isometry group, standard operator algebra, extensions of $C^*$-algebras,
Toeplitz algebra, Laurent algebra
\endkeywords
\thanks
This paper was written when the author, holding a scholarship
of the Volkswagen-Stiftung of the Konferenz der Deutschen Akademien der
Wissenschaften, was a
visitor at the University of Paderborn, Germany. He is very grateful to
Prof. K.-H. Indlekofer for his kind hospitality.
This research was partially supported also by the Hungarian National
Foundation for Scientific Research (OTKA), Grant No. T--016846
F--019322.
\endthanks
\address
Institute of Mathematics, Lajos Kossuth University, 4010
Debrecen, P.O.Box 12, Hungary
\endaddress
\email
{\tt molnarl\@math.klte.hu}
\endemail
\date
April 25, 1997
\enddate
\abstract
In this paper we answer a question raised in \cite{Mol}
by giving an example of a proper standard $C^*$-subalgebra
of $\B(\H)$ whose automorphism and isometry groups are topologically
reflexive. After this we prove that in the case of extensions of the
$C^*$-algebra $\C(\H)$ of all compact operators by separable
commutative $C^*$-algebras, these groups are algebraically
reflexive. Concerning the most well-known extension of $\C(\H)$ by
$C(\Bbb T)$ (the algebra of all continuous complex valued functions
on the perimeter of the unit disc)
we show that the automorphism and isometry groups
are topologically nonreflexive.
\endabstract
\endtopmatter

\document
\head
1. Introduction and Statements of The Results
\endhead

Let $\X$ be a Banach space and let $\B(\X)$ denote the
algebra of all bounded linear operators on $\X$.
A subset $\E \subset \B(\X)$ is called topologically [algebraically]
reflexive if for every $T\in \B(\X)$, the condition that $Tx\in
\overline{\E x}$(=the norm-closure of $\E x$) [$Tx\in \E x$] holds true
for every $x\in \X$ implies $T\in \E$ (cf. \cite{LoSu}).
Roughly speaking, reflexivity means that the elements of $\E$
are, in some sense, completely determined by their local actions.
The concept of reflexive subspaces has been proved very
useful in the analysis of operator algebras (see, for example,
\cite{Lar} and the references therein).

Reflexivity problems concerning sets of linear transformations on
operator algebras were first studied by Kadison and Larson and Sourour.
In \cite{Kad}, \cite{LaSu} the problem of algebraic reflexivity of the
linear space of all derivations on a von Neumann algebra,
respectively on $\B(\X)$ was discussed and solved.
As for topological reflexivity, Shul'man \cite{Shu}
proved that the derivation algebra of any $C^*$-algebra is
topologically reflexive.
In \cite{BrSe} Bre\v sar and \v Semrl showed that, for a separable
Hilbert space $\H$, the set of all automorphisms of
$\B(\H)$ is algebraically reflexive.
Here we emphasize that in the present paper automorphism means a
merely multiplicative
linear bijection, so the *-preserving property is not supposed.
In our papers \cite{BaMo, Mol} we proved that the group of
all automorphisms as well as the group of all surjective isometries
of $\B(\H)$ are topologically reflexive.
Since the topological reflexivity of these groups seemed to be a
quite exceptional phenomenon,
in \cite{Mol} we raised the question of the existence of a proper
$C^*$-subalgebra of $\B(\H)$ containing all compact operators whose
automorphism and isometry groups are
topologically reflexive. In this present paper we show that the answer
to this question is affirmative. In fact, we present a quite
simple algebra as an example. Afterwards, we
discuss the reflexivity of the groups in question in the case of
extensions of the $C^*$-algebra $\C(\H)$ of compact operators which are
of
fundamental importance in operator theory. These are the extensions of
$\C(\H)$ by separable commutative $C^*$-algebras which are the main
objects of the celebrated Brown-Douglas-Fillmore theory
\cite{Dav, Chapter IX}. The most well-known examples of these extensions
are the so-called Toeplitz algebra and Laurent algebra
both of which are extensions of $\C(\H)$ by $C(\Bbb T)$.
We show that the automorphism and isometry groups of these algebras are
algebraically reflexive but they fail to be topologically reflexive.

In what follows let $\H, \K$ be infinite dimensional separable complex
Hilbert spaces.
Denote by $\F$
the ideal of all finite rank operators in $\B(\H)$.
A subalgebra $\A$ of $\B(\H)$ is called standard if it contains
$\F$.
Let $\{ \H_i\}$ be a fixed (finite or infinite) sequence of pairwise
orthogonal closed subspaces of $\H$ which generate $\H$. Let $\B(\{
\H_i\})$ denote the $C^*$-subalgebra of $\B(\H)$ consisting of all
operators $A\in \B(\H)$ for which $A(\H_i)\subset \H_i$ for every $i$.

Our first theorem answers the open problem raised at the end of our
former paper \cite{Mol} concerning the existence of a proper standard
$C^*$-subalgebra of $\B(\H)$ whose automorphism and isometry groups are
topologically reflexive.

\proclaim{Theorem 1.1}
Let the number of the subspaces $\{ \H_i\}$ be finite.
Then the automorphism group and the isometry group of the
$C^*$-algebra $\C(\H)+\B(\{ \H_i\})$ are topologically reflexive.
\endproclaim

The remaining results of the paper treat our reflexivity problem in the
case of extensions of $\C(H)$ by
separable commutative $C^*$-algebras. Let $X$ be a compact metric space.
Denote by $C(X)$ the $C^*$-algebra of all continuous complex valued
functions on $X$. Recall that these are exactly the
separable unital commutative $C^*$-algebras. By an
extension of $\C(\H)$ by $C(X)$ we mean a standard $C^*$-algebra $\A$
together with a surjective *-homomorphism $\phi :\A \to C(X)$ whose
kernel $\ker \phi=\C(\H)$.
The most well-known such extensions are concerned with the case $X=\Bbb
T$. The building blocks of these extensions are the
Toeplitz algebra and the Laurent algebra.
The Toeplitz algebra $\eusm{T}$ is defined by
$$
\eusm{T}=\{ K+T_f \: K\in \C(\Bbb L^2(\Bbb T)), f\in C(\Bbb T)\},
$$
where $\Bbb T$ is the perimeter of the unit disc with normalized
Lebesgue measure and $T_f$ denotes the Toeplitz operator corresponding
to the continuous symbol $f$ (see \cite{Cob1-2}, \cite{Dav, Sections
V.1 and V.2}).
This can also be realized as the $C^*$-algebra generated by a unilateral
shift $S$. The Laurent algebra $\eusm{L}$ is defined by
$$
\eusm{L}=\{ K+M_f \: K\in \C(\Bbb L^2(\Bbb T)),
f\in C(\Bbb T)\},
$$
where $M_f$ denotes the operator of multiplication by the function $f$.
This can also be realized as the $C^*$-algebra generated by the compact
operators and a bilateral shift. Concerning the reflexivity of the
automorphism
and isometry groups of these algebras we have the following results.

\proclaim{Theorem 1.2}
Let $X$ be a compact metric space. The automorphism and isometry
groups of any extension of $\C(\H)$ by $C(X)$ are algebraically
reflexive.
\endproclaim

\proclaim{Theorem 1.3}
Let $\A$ be any extension of $\C(\H)$ by $C(\Bbb T)$.
Then there is a sequence of *-automorphisms of $\A$
which converges pointwise to a nonsurjective *-endomorphism.
Therefore, the automorphism group as well as the isometry group
of $\A$ are topologically nonreflexive.
\endproclaim

The paper is concluded by some remarks and open problems.

\head
2. Proofs
\endhead

Let us begin with some observations, remarks and notation which we make
use in the proofs. A continuous linear map $\Cal J$ between
normed algebras $\Cal A$ and $\Cal B$ is called a Jordan homomorphism if
$$
\Cal J(A)^2=\Cal J(A^2) \qquad (A\in \Cal A).
$$
Observe that linearizing the previous equality, i.e. replacing $A$ by
$A+B$ we can deduce that $\Cal J$ satisfies
$$
\Cal J(AB+BA)=\Cal J(A)\Cal J(B)+\Cal J(B)\Cal J(A) \qquad (A,B \in \Cal
A).
$$

Our main objectives are the standard $C^*$-algebras. The structures of
all Jordan automorphisms, automorphisms, antiautomorphism (i.e. linear
bijections reversing the order of multiplication) as well
as surjective isometries of these
algebras are well-known
and easy to describe as we see in the following proposition.

\proclaim{Proposition 2.1}
Let $\Cal A\subset \B(\H)$ be a standard $C^*$-algebra.
Then every Jordan automorphism of $\A$ is either an automorphism or an
antiautomorphism. In the first case we have an invertible bounded linear
operator $T$ on $\H$ such that $\P$ is of the form
$$
\P(A)=TAT^{-1} \qquad (A\in \Cal A).
$$
In the second case we have an invertible bounded linear operator $S$ on
$\H$ such that $\P$ is of the form
$$
\P(A)=SA^{tr}S^{-1} \qquad (A\in \Cal A)
$$
where ${}^{tr}$ denotes the transpose with respect to an arbitrary but
fixed complete
orthonormal sequence in $\H$. This latter assertion is equivalent to
saying that there is an invertible bounded conjugate-linear operator
$S'$ on $\H$ such that
$$
\P(A)=S'A^{*}{S'}^{-1} \qquad (A\in \Cal A).
$$
If $\Cal A$ contains $I$ and $\S:\Cal A\to \Cal A$ is a surjective
linear isometry, then there
are unitary operators $U,V$ on $\H$ such that $\P$ is either of the form
$$
\S(A) =UAV\qquad (A\in \Cal A)
$$
or of the form
$$
\S(A) =UA^{tr}V\qquad (A\in \Cal A).
$$
\endproclaim

\demo{Proof}
It is a well-known theorem of Herstein \cite{Her} that every Jordan
homomorphism onto a prime algebra is either a homomorphism or an
antihomomorphism. Since every standard operator algebra is prime, we
have the first assertion.
It is a classical theorem of Kadison \cite{KaRi, 7.6.17, 7.6.18}
that every surjective isometry of a unital $C^*$-algebra is a Jordan
*-automorphism (i.e. a Jordan automorphism preserving the *-operation)
multiplied by a fixed unitary element. Now, our statement follows from
folk results on the forms of automorphisms,
antiautomorphisms, *-automorphisms and
*-antiautomorphisms of standard operator
algebras (cf. \cite{Che}).
\enddemo

Let $\P:\Cal A\to \Cal B$ be an approximately local
homomorphism, i.e. a continuous linear map such that for every $A\in
\Cal A$ there is a sequence $(\P_n)$ of homomorphisms
(depending on $A$) for which $\P(A)=\lim_n \P_n(A)$. We are
interested in the question when it follows that $\P$ is a
Jordan homomorphism. The easy proposition below based on a
well-known computation will be of some help in what follows.
Observe that every approximately local homomorphism sends
projections to idempotents.

\proclaim{Proposition 2.2}
Let $\Cal A, \Cal B\subset \B(\H)$ be closed *-subalgebras and suppose
that for every self-adjoint element $A$ of $\A$, the spectral measure of
any Borel subset of $\sigma (A)$ bounded away 0 belongs to
$\A$. If $\P:\Cal A\to \Cal B$ is a continuous linear map
which sends projections to idempotents, then $\P$
is a Jordan homomorphism.
\endproclaim

\demo{Proof}
Let $P,Q\in \Cal A$ be mutually orthogonal projections.
Then $\P(P)+\P(Q)$ is an idempotent and since it is the sum of two
idempotents, we have $\P(P)\P(Q)=\P(Q)\P(P)=0$. If
$\l_1,\dots ,\l_n$ are real numbers and $P_1,\dots ,P_n\in \Cal A$ are
mutually orthogonal projections, then we infer
$$
(\P(\sum_{k=1}^n \l_k P_k))^2=
(\sum_{k=1}^n \l_k \P(P_k))^2=
\sum_{k=1}^n \l_k^2 \P(P_k)=
\P((\sum_{k=1}^n \l_k P_k)^2).
$$
By the spectral theorem and the continuity of $\P$ this implies that
$\P(A)^2=\P(A^2)$ holds true for every self-adjoint element $A\in \Cal
A$. Linearizing this equality,
we immediately get $\P(AB+BA)=\P(A)\P(B)+\P(B)\P(A)$
for any self-adjoint $A,B\in \Cal A$.
Finally, if $T\in \Cal A$ is arbitrary, then it can be written in the
form $T=A+iB$ with self-adjoint $A,B\in \Cal A$ and
the previous equalities result in the desired $\P(T)^2=\P(T^2)$.
\enddemo

Now, considering the statement of Theorem 1.1,
if the algebra $\C(\H)+\B(\{ \H_i\})$ had the property concerning
spectral measures appearing in the formulation of Proposition 2.2,
then we could infer that every approximately local
automorphism of $\C(\H)+\B(\{ \H_i\})$ is a Jordan homomorphism.
Unfortunately, this algebra fails to have that property.
Indeed,
let the number of the subspaces $\{ \H_i\}$ be at least 2. Consider
an arbitrary infinite dimensional projection $P\in \B(\H)$. Clearly,
$P$ can be written as $P=E([1,\infty[\cap \sigma(I+K))$ for some
positive compact operator $K$.
Hence, if the algebra $\C(\H)+\B(\{ \H_i\})$ had the property
in Proposition 2.2, then
we would obtain that every projection in $\B(\H)$ belongs to
$\C(\H)+\B(\{ \H_i\})$. Apparently, this implies $\C(\H)+\B(\{
\H_i\})=\B(\H)$ which is a contradiction.

The second thought which one might have is that
it is possible to get Theorem 1.1 directly from
\cite{Mol, Theorem 2 and Theorem 3}. Namely, it may be guessed that the
almost only thing that we have to do is to
verify that every approximately local automorphism of the $C^*$-algebra
$\C(\H)+\B(\{ \H_i\})$ is, by restriction, an approximately local
automorphism of $\B(\{ \H_i\})$ which is the direct sum of $\B(\H_i)$'s
and we could try to apply our result \cite{Mol, Theorem 2} on the
topological reflexivity of the automorphism group of $\B(\H)$.
Once again, the starting point of this argument is false.
In fact, it is easy to give an example of an automorphism of
$\C (\H \oplus \H)+ \B(\{\H, \H \})$ whose restriction to
$\B(\{\H, \H \})$ is not an automorphism. For instance,
let $0\neq P$ be a finite rank projection on $\H$ and let $S,T\in
\B(\H)$
be such that $\ker S=\ker T=\rng P$, $\rng S=\rng T=\rng (I-P)$ and
$ST=TS=I-P$. Using elementary computations, one can verify
that the map
$$
\left[
\matrix
A & K\\
C & B
\endmatrix
\right ]
\longmapsto
\left[
\matrix
P & T\\
S & P
\endmatrix
\right ]
\left[
\matrix
A & K\\
C & B
\endmatrix
\right ]
\left[
\matrix
P & T\\
S & P
\endmatrix
\right ]
$$
is an automorphism of the algebra $\C (\H \oplus \H)+ \B(\{\H, \H \})$
which does not leave $\B(\{ \H, \H \})$ invariant.

Now, we turn to the proof our first theorem which we reach via a series
of auxiliary statements.

\proclaim{Lemma 2.3}
Let $Q_n$ be a bounded sequence of idempotents in $\B(\H)$ such that
$\rng Q_n \subset \rng Q_{n+1}$ $(n\in \N)$. Then $(Q_n)$
converges strongly to an idempotent $Q\in \B(\H)$.
\endproclaim

\demo{Proof}
Elementary functional analysis.
\enddemo

\proclaim{Proposition 2.4}
Let $\A\subset \B(\H)$ be a standard $C^*$-algebra which
is linearly generated (in the norm topology) by
the set of its projections and suppose that for every closed nontrivial
ideal
$\I$ of $\A$, the quotient algebra $\A/ \I$ contains uncountably many
pairwise orthogonal projections.
Let $\P:\A\to \B(\K)$ be a Jordan homomorphism. If $(P_n)$ is a
maximal family of rank-one projections in $\B(\H)$, then
the idempotent $E=\sum_n \P(P_n)$ is well-defined (we mean that it does
not depend on the particular choice of $(P_n)$), $E$ commutes with
the range of $\P$ and we have $\P(.)=\P(.)E$.
\endproclaim

\demo{Proof}
The assertions that $E$ is well-defined and it commutes with
the range of $\P$ follow directly from the proof of \cite{Mol, Lemma 2}.
As for the remaining statement $\P(.)=\P(.)E$, observe that
the map
$$
\Psi: A\longmapsto \P(A)(I-E)
$$
is a Jordan homomorphism and it is easy to see that $\Psi$
vanishes on every finite-rank projection.
The kernel $\I$ of $\S$ is a closed Jordan ideal of
$\A$. It is well-known that every closed Jordan ideal in a $C^*$-algebra
is an associative ideal as well \cite{CiYo}. Therefore, $\I$ is a closed
associative ideal in $\A$. If $\I\neq \A$, then by our assumption
on $\A$ it follows that the range of $\S$ contains an uncountable family
of pairwise orthogonal nonzero idempotents. Since this contradicts the
separability of $\K$, we have $\S=0$. This gives us that
$\P(.)=\P(.)E$.
\enddemo

\proclaim{Corollary 2.5}
Let $\A$ be as in Proposition 2.4 above. If
$\P, \P':\A \to \B(\K)$ are
Jordan homomorphisms which coincide on $\F$, then $\P=\P'$.
\endproclaim

\demo{Proof}
Let $(P_n)$ be a maximal family of pairwise orthogonal
rank-one projections in $\B(\H)$. Let $A\in \A$ be arbitrary.
By Proposition 2.4 we infer
$$
\gather
2\P(A)=\sum_n (\P(A)\P(P_n) +\P(P_n)\P(A))=
\sum_n \P(AP_n +P_nA)=\\
\sum_n \P'(AP_n +P_nA)=
\sum_n (\P'(A)\P'(P_n) +\P'(P_n)\P'(A))=2\P'(A).
\endgather
$$
\enddemo

\proclaim{Proposition 2.6}
Let $\P :\C(\H) \to \C(\K)$ be a Jordan homomorphism.
Then the second adjoint $\P^{**}$ of $\P$ defines a weak*-continuous
Jordan homomorphism from $\B(\H)$ to $\B(\K)$ which extends
$\P$.
\endproclaim

\demo{Proof}
It is well-known that the dual space of $\C(\H)$ is the Banach algebra
$\eusm{T}(\H)$ of all trace-class operators on $\H$ and the dual space
of $\eusm{T}(\H)$ is $\B(\H)$. The dualities in
question are given by the pair
$$
\la A,B\ra =\tr AB
$$
where $A\in \C(\H), B\in \eusm{T}(\H)$, respectively $A\in
\eusm{T}(\H), B\in \B(\H)$. Here, $\tr$ denotes the usual
trace-functional.

Now, if $K\in \C(\H)$ and $T\in \eusm{T}(\K)$, then we
compute
$$
\tr \P^{**}(K)T=\tr K\P^{*}(T)=\tr \P(K)T.
$$
This apparently gives us that $\P^{**}$ is an extension of $\P$.
Let $P\in \B(\H)$ be an arbitrary projection and let $(P_n)$ be a
monoton increasing sequence of finite rank projections which
converges strongly to $P$. We then have $\tr P_n T \to \tr PT$ for every
trace-class operator $T$ and we infer
$$
\tr \P(P_n)T=\tr P_n \P^*(T) \longrightarrow \tr P\P^*(T)= \tr
\P^{**}(P)T.
$$
This implies that $\P(P_n)$ converges weakly to $\P^{**}(P)$. On the
other hand, by Lemma 2.3 it follows that $\P(P_n)$ converges
strongly to an idempotent. Hence, $\P^{**}(P)$ is an idempotent whenever
$P$
is a projection. Using Proposition 2.2 we obtain that $\P^{**}$ is a
Jordan homomorphism.
\enddemo

\proclaim{Proposition 2.7}
Let $\A$ be as
in Proposition 2.4. If, in addition, $\A$ contains $I$, and $\P:\A \to
\B(\K)$ is a unital Jordan homomorphism which
preserves the rank-one operators, then there is an invertible bounded
linear operator $T:\H \to \K$ so that $\P$ is either of the form
$$
\P(A)=TAT^{-1} \qquad (A\in \A)
$$
or of the form
$$
\P(A)=TA^{tr}T^{-1} \qquad (A\in \A).
$$
\endproclaim

\demo{Proof}
Let $\S$ be the restriction of $\P$ onto $\C(\H)$. Clearly, $\rng
\S\subset \C(\H)$. By Proposition 2.7, $\S^{**}$ is a weak*-continuous
Jordan homomorphism which preserves the rank-one
operators. We now apply a result of Hou \cite{Hou, Theorem 1.2} on the
form of linear maps sending rank-one operators to operators with rank at
most one. It says that either there are continuous linear operators
$T:\H\to \K$ and $S:\K \to \H$ so that $\S^{**}$ is of the form
$$
\S^{**}(A)=TAS \qquad (A\in \B(\H))
$$
or there are bounded conjugate-linear operators $T':\H \to \K$ and
$S':\K \to \H$ so that $\P$ is of the form
$$
\S^{**}(A)=T'A^*S' \qquad (A\in \B(\H)).
$$
In fact, Hou's theorem was formulated for weak-continuous maps but
an inspection of the proof shows that
this condition was used only to prove the continuity of $T,S$ and
to show that if the above formula is valid on $\F$, then it holds
true
on $\B(\H)$ as well. Obviously, in both places weak*-continuity can play
the same role. Going further in our proof, let us suppose that
$\S^{**}$ is of the first form.
By Corollary 2.5 and Proposition 2.6 we have $\S^{**}_{|\A}=\P$.
Since $\S^{**}$ is a Jordan homomorphism, it preserves the
idempotents. Therefore, if $x\ot y$ is such that $\la x,y\ra =1$, then
we
have $\la Tx,S^* y\ra=1$. Clearly, it implies $\la STx,y\ra =\la x,y\ra$
$(x,y\in \H)$ which gives us that $ST=I$. On the other hand,
$\P$ is unital and hence we infer $TS=I$. Consequently,
$S=T^{-1}$.
If $\P^{**}$ is of the second form above, one can follow the same
argument.
\enddemo

After this preparation we now are in a position to prove our first
theorem.

\demo{Proof of Theorem 1.1}
So, let the number of the subspaces $\{ \H_i\}$ be finite. Suppose that
$\P$ is a continuous linear map which is an approximately local
automorphism of the $C^*$-algebra $\A=\C(\H)+\B(\{ \H_i\})$. We
first show that $\P$ is a Jordan homomorphism. Clearly, the
restrictions $\P_{|\C(\H)}$ and $\P_{|\B(\{ \H_i\})}$ send idempotents
to idempotents.
Proposition 2.2 guarantees that they are Jordan homomorphisms.
For every $i$ we define a linear map $\P_i :\B(\H_i) \to \B(\K)$ by
$$
\P_i (A)=\P(\hat A)
$$
where $\hat A$ is the element of $\B(\{ \H_i\})$ which coincide with
$A$ on $\H_i$ and 0 on $\H_i^\perp$.
Furthermore, let $\S=\P_{|\C(\H)}^{**}$. Here we note that, by
the local property of $\P$, it sends rank-one operators to operators
with rank not greater than one. In fact, $\P$ is a
rank-one preserver. To see this, observe that
if we suppose on the contrary that $\P$ sends a rank-one operator to
0, then examining the kernel of $\P$ we would get that $\P$ vanishes
on $\C(\H)$. Now, if $\H_i$ is any infinite dimensional subspace in
our collection $\{ \H_i\}$, we can infer that $\P$ vanishes on
$\B(\H_i)$. Obviously,
this results in $\P=0$ which is a contradiction. Consequently,
$\P$ is a rank-one preserver. Let us define
$$
\S_i(A)=\S(\hat A) \qquad (A\in \B(\H_i)).
$$
The maps $\P_i$, $\S_i$ are Jordan homomorphisms which coincide on
$\C(\H_i)$. Therefore, by Corollary 2.5 it follows that they are equal
which, after summation, gives us that
$$
\P(A)= \S(A)
$$
for every $A\in \B(\{ \H_i\})$. Since $\S$ is a
Jordan homomorphism which extends $\P_{|\C(\H)}$, we compute
$$
\gather
\P(A)\P(K)+\P(K)\P(A)=
\S(A)\S(K)+\S(K)\S(A)=\\
\S(AK+KA)=
\P(AK+ KA)
\endgather
$$
for every $A\in \B(\{ \H_i\}), K\in \C(\H)$.
The fact that $\P$ is a Jordan homomorphism now follows from the
equality
$$
\gather
(\P(K+A))^2=\P(A)^2+\P(A)\P(K)+\P(K)\P(A)+\P(K)^2=\\
\P(A^2+AK+KA+K^2)= \P((K+A)^2).
\endgather
$$
Next, since $\P$ is unital and preserves the rank-one operators, from
Proposition 2.7 it follows that
there is an invertible operator $T\in \B(\H)$ such that $\P$ is either
of the form
$$
\P(A)=TAT^{-1} \qquad (A\in \A)
$$
or of the form
$$
\P(A)=TA^{tr}T^{-1} \qquad (A\in \A).
$$
Suppose that $\P$ is of this latter form. Let $\H_i$ be an infinite
dimensional subspace from our collection $\{ \H_i\}$. Pick an operator
$U\in \A$ which is unilateral shift on $\H_i$
and the identity on $\H_i^\perp$. Obviously, $U$ has a left inverse in
$\A$ but it does not have a right one. Clearly, the same must hold true
for the image of $U$ under any automorphism of $\A$. Let $(\P_n)$ be a
sequence of automorphisms of $\A$ for which $\P(U)=\lim_n \P_n(U)$.
Since the set of all elements which have right inverse in $\A$ is
open, we deduce that $\P(U)$ has no right inverse. On the other hand, we
can compute
$$
\P(U)\P(U^*)=\P(U^*U)=\P(I)=I.
$$
Thus, we have arrived at a contradiction and, consequently, it follows
that $\P(A)=TAT^{-1}$ for every $A\in \A$.

We know that $\P:\A \to \A$ and hence $TAT^{-1}\in \A$ holds true for
every $A\in \A$. We claim that this implies that $T^{-1}AT\in \A$
$(A\in \A)$ which then will give us the surjectivity of $\P$.
Consider the matrix representations
of the elements of $\A$ corresponding to the subspaces $\{ \H_i\}$. Let
$T=[T_{ij}]$ and $T^{-1}=S=[S_{ij}]$. Let the index $i_0$ be fixed for
a moment and pick any operator $A_{i_0}\in \B(\H_{i_0})$. By
$T\A T^{-1}\subset \A$ we obtain that the off-diagonal elements
of the matrix $[T_{ii_0}A_{i_0}S_{i_0j}]$ are all compact
operators. So, for any $i\neq j$ we have
$T_{ii_0}\B(\H_{i_0})S_{i_0j}\subset \C(\H_j, \H_i)$.
Using, for example, the characterization of compact operators as those
bounded linear operators whose range does not contain any infinite
dimensional closed subspace, it is easy to see that we necessarily have
that either $T_{ii_0}$ or $S_{i_0j}$ must be compact.
Now, let us remove those rows and colums from the matrices of $T$ and
$S$ which correspond to finite dimensional subspaces but hold on
the numbering of the entries.
Denote the
matrices obtained in this way by $\tilde T$ and $\tilde S$,
respectively. Obviously,
we still have the property that, considering the $i$th column of $\tilde
T$ and the $i$th row of $\tilde S$, from any pair of entries sitting in
different positions, one of them is compact. We show that in every
row and column of $\tilde T$ there is exactly one non-compact entry and
the same holds true for $\tilde S$. To see this, consider the $i$th
column of $\tilde
T$. If every entry of it is compact, then by $ST=I$ it follows that the
identity on $\H_i$ is compact which implies that $\H_i$ is finite
dimensional and this is a
contradiction. Next, suppose that there are two non-compact entries
in the column in question. Then it easily follows that in the
$i$th row of
$\tilde S$ consists of compact entries. Using $ST=I$ just as above, we
arrive at a contradiction again. Therefore, there is exactly one
non-compact
entry in every column of $\tilde T$. Suppose that there is a row in
$\tilde T$ which contains two non-compact elements. Then we
necessarily have another row
of $\tilde T$ whose entries are all compact. But by $TS=I$
this is untenable. Hence, we have proved that
every row and column of $\tilde T$ contains exactly one non-compact
entry. Clearly,
this implies that $\tilde S$ has the same property. In fact, there is a
non-compact element in position $ij$ in $\tilde T$ if and only if there
is
a non-compact element in position $ji$ in $\tilde S$. Now, it is
apparent that the off-diagonal elements in
$[S_{ii_0}A_{i_0}T_{i_0j}]$ are all compact.
This gives us the desired inclusion $S\A T\subset \A$ and we obtain
the topological reflexivity of the automorphism group.

Let us now prove the topological reflexivity of the isometry group.
By Proposition 2.1 every surjective isometry of $\A$ preserves the
unitary group. Plainly, if $\P$ is a
continuous linear map which is an approximately local surjective
isometry, then $\P$ has the same preserver property. But the
structure of unitary group preservers of $C^*$-algebras is well-known.
In fact, \cite{RuDy, Corollary} gives us that there is a unital Jordan
*-homomorphism $\S$ on $\A$ and a unitary element $U\in \A$ so that
$\P(A)=U\S(A)$ $(A\in \A)$. Obviously, we may suppose that $U=I$.
From the form of Jordan automorphisms of standard operator algebras it
follows that $\P$ preserves the rank-one operators.
Therefore, by Proposition 2.7 we infer that there is an
invertible operator (in fact, a unitary one in the case of
Jordan *-homomorphisms) $T$ such that $\P$ is either of the form
$$
\P(A)=TAT^{-1} \qquad (A\in \A)
$$
or of the form
$$
\P(A)=TA^{tr}T^{-1} \qquad (A\in \A).
$$
This latter form can be rewritten as $\P(A)=T'A^*{T'}^{-1}$
with some invertible bounded conjugate-linear operator $T'$.
The proof can be completed as in the case of the automorphism group.
\enddemo

We now turn to the proofs of our results on the extensions of $\C(\H)$.

\proclaim{Proposition 2.8}
Let $X$ be a I. countable compact Hausdorff space. The
automorphism and isometry groups of $C(X)$ are algebraically reflexive.
\endproclaim

\demo{Proof}
The algebraic reflexivity of the isometry group of $C(X)$
was proved in \cite{MoZa, Theorem 2.2}.
As for the case of automorphisms, we recall that the automorphisms
of $C(X)$ are
of the form $f\mapsto f\circ \varphi$, while the surjective isometries
are of the form $f\mapsto \tau f\circ \varphi$, where $\varphi:X\to X$
is a homeomorphism and $\tau$ is a continuous complex valued function on
$X$ with absolute value 1.
The algebraic reflexivity of the automorphism group now easily follows
from that of the isometry group.
\enddemo

\proclaim{Proposition 2.9}
The automorphism and isometry groups of $\C(\H)$ are algebraically
reflexive.
\endproclaim

\demo{Proof}
Let $\P:\C(\H)\to \C(\H)$ be a continuous linear map which is a local
automorphism of $\C(\H)$. By the form of the automorphisms of $\C(\H)$
(see Proposition 2.1), it is obvious that $\P$ preserves the rank-one
operators.
Using Hou's result again with additional remarks very similar to those
ones which appeared in the proof of Proposition 2.7,
it follows that either there are bounded linear operators $T,S:\H \to
\H$ so that $\P$ is of the form
$$
\P(A)=TAS \qquad (A\in \C(\H))
$$
or there are bounded conjugate-linear operators $T',S':\H \to \H$ so
that $\P$ is of the form
$$
\P(A)= T'A^*S' \qquad (A\in \C(\H)).
$$
Since $\P$ is a Jordan homomorphism (see Proposition 2.2), following the
same argument as in the proof of Proposition 2.7, it
is easy to conclude that $ST=I$, respectively $S'T'=I$ hold true.
We show that the appearence of the second form of $\P$ can
be excluded. In fact, if $A$ is a compact
operator such that $A$ is injective but $A^*$ is not, then by the local
property of $\P$ it follows that the same must be valid for $\P(A)$ as
well. But if $T'A^*S'$ is injective, then using the surjectivity of $S'$
what we get from $S'T'=I$, we obtain that $A^*$ is injective. Since
this is a contradiction, we infer that our $\P$ is of the first form.
Considering $ST=I$ again and referring to the fact that $TAS=\P(A)$
is injective whenever $A$ is so, it is apparent that $S=T^{-1}$. This
completes the proof in the case of automorphisms.

The algebraic reflexivity of the isometry group of $\C(\H)$ was
proved in \cite{MoZa, Theorem 1.6}.
\enddemo

\demo{Proof of Theorem 1.2}
Let $\A$ be an extension of $\C(\H)$ by $C(X)$ and let
$\P :\A \to \A$ be a local automorphism. By the form  of the
automorphisms of standard operator algebras, it follows that the
restriction of any automorphism of $\A$ onto $\C(\H)$ is an automorphism
of $\C(\H)$. Therefore, considering quotients and using
the identification of $\A/\C(\H)$ and $C(X)$, we get that every
automorphism of $\A$ gives rise to an automorphism of $C(X)$.
Let $\phi (f)=\overline{\P(
A)}$, where $A\in \A$ is such that $\overline{A}=f$. Since $\phi$ is a
local automorphism of $C(X)$, by Proposition 2.8 it follows that $\phi$
is surjective, that is
$\{ \overline{\P(A)} \: A\in \A\}=\A/\C(\H)$. As a consequence, we
obtain that to every operator $B\in \A$ there corresponds an operator
$A\in \A$ so that $\P(A)-B\in \C(\H)$. Since the restriction of $\P$
onto
$\C(\H)$ is a local automorphism of $\C(\H)$, by Proposition 2.9 it is
an automorphism and hence $\P(A)-B=\P(K)$ holds true for
some
compact operator $K$. Clearly, $B=\P(A-K)$ and we obtain the
surjectivity of $\P$.
Now, we can apply a very nice result of Aupetit and Mouton. Since our
mapping $\P$ is a local automorphism, it preserves the
invertible operators in both directions. This shows that $\P$ is a
spectrum-preserving map onto the primitive Banach algebra $\A$ which
contains a minimal ideal. By \cite{AuMo, Corollary 3.4} it follows that
$\P$ is either an automorphism or an antiautomorphism.
Using the same argument as in the proof of Proposition 2.9, it quickly
follows that $\P$ is an automorphism
proving the algebraic reflexivity of the automorphism group of $\A$.

Now, let $\P$ be a local surjective isometry of $\A$. Since $\P$ is
automatically an isometry, only the surjectivity needs proof.
But taking the form of surjective isometries of standard operator
algebras into consideration, this can be derived just as in the case
of local automorphisms.
\enddemo

\demo{Proof of Theorem 1.3}
Let us first consider the case of the Laurent algebra.
It is easy to see that $\eusm{L}$ can be viewed also in the
following way
$$
\eusm{L}=\{ K+M_f\: K\in \C(\Bbb L^2[0,2\pi]), f\in C[0,2\pi],
f(0)=f(2\pi)\}.
$$
Define a sequence $(\vp_n)$ of homeomorphisms of $[0,2\pi]$ in such a
way that the following conditions be fulfilled:

\itemitem{(i)} $\vp_n(0)=0$, $\vp_n(2\pi)=2\pi$,

\itemitem{(ii)} $\vp_n$ is continuously differentiable and its
derivative vanishes nowhere,

\itemitem{(iii)} $(\vp_n)$ converges uniformly to the function $\vp$
defined by $\vp(x)=2x$ $(x\in [0,\pi])$, $\vp(x)=2\pi$ $(x\in ]\pi,
2\pi])$,

\itemitem{(iv)} $\vp'_n(x)\to 2$ $(x\in [0,\pi[)$ and
$\vp'_n(x)\to 0$ $(x\in ]\pi,2\pi])$.

\noindent
Define operators $V_n, V\in \B(\Bbb L^2[0,2\pi])$ by
$$
\split
V_n g&=\sqrt{\vp_n'} (g\circ \vp_n)\\
Vg   &=\cases \sqrt 2 g(2x) \, \text{ if } \, x\in [0,\pi]\\
                0 \phantom{aig(2x)} \, \text{ if } \, x\in ]\pi, 2\pi]
\endcases
\endsplit
$$
$(g\in \Bbb L^2[0,2\pi])$.
It is elementary to verify that $V_n$ is unitary $(n\in \N)$, $V$ is a
nonsurjective isometry and $(V_n)$ converges strongly to $V$.
An easy computation shows that $M_f =V_n^*M_{f\circ \vp_n}V_n$
whenever $f\in C[0,2\pi]$, $f(0)=f(2\pi)$. Then we have $M_{f\circ
\vp_n}=V_nM_fV_n^*$ and it is trivial to check that the formula
$$
\P_n(K+M_f)=V_n(K+M_f)V_n^* =V_nKV_n^*+M_{f\circ \vp_n}
$$
$(K\in \C(\Bbb L^2[0,2\pi]), f\in C[0,2\pi], f(0)=f(2\pi))$
defines a sequence of *-auto\-mor\-phisms of the Laurent algebra. Let us
define $\P:\eusm{L}\to \eusm{L}$ by
$$
\P(K+M_f)=VKV^*+M_{f\circ \vp}
$$
for every $K\in \C(\Bbb L^2[0,2\pi]), f\in C[0,2\pi], f(0)=f(2\pi)$.
Since $(V_n)$ converges strongly to $V$, we obtain that
$V_n FV_n^* \to VFV^*$ holds true in the operator norm topology for
every finite
rank operator and hence, by Banach-Steinhaus theorem, for every compact
operator as well.
By the uniform convergence $\vp_n \to \vp$, we deduce that
$M_{f\circ \vp_n}\to M_{f\circ \vp}$.
Therefore, we have $\P_n(A)\to \P(A)$ $(A\in \eusm{L})$.
Finally, since $V$ is a proper isometry and every nonzero multiplication
operator is
non-compact, it is obvious that the range of $\P$ does not contain every
compact operator.

Now, we turn to the case of the Toeplitz algebra.
Let $S=\sum_{n=1}^\infty e_{n+1}\ot e_n$ be the generating
unilateral shift,
where $(e_n)$ is a complete orthonormal sequence in $\H$. Let $S'=e_1\ot
e_1+ \sum_{n=2}^\infty e_{n+1}\ot e_n$. Since $S'$ is the direct sum of
a (one-dimensional) unitary and a unilateral shift, by \cite{Hal,
Lemma 5} there is a sequence $(U_n)$ of unitaries on $\H$ such that
$U_n SU_n^* \to
S'$ in the operator norm and $U_nSU_n^*-S'$ is compact for every $n$.
Since $S'$ is a finite-rank perturbation of $S$, we have $S'\in \eusm{T}$
and then $U_nSU_n^*\in \eusm{T}$.
Since $U_n SU_n^* $ is a unilateral shift, the
$C^*$-algebra $\A(U_nSU_n^*)$ generated by $U_nSU_n^*$ contains every
compact operator and by $U_nSU_n^*-S\in \C(\H)$ it follows that $S\in
\A(U_nSU_n^*)$. Therefore, we have $\A(U_nSU_n^*)=\eusm{T}$.
Since $S'\in \eusm{T}$, we obtain
$\A(S')\subset \eusm{T}$ but the converse
inclusion does not hold true due to the fact that $S'$ and
hence every operator in $\A(S')$ has $e_1$ as an eigenvector.
Let us define $\P(S)=S'$ and $\P_n(S)=U_nSU_n^*$.
By Coburn's theorem \cite{Dav, Theorem V.2.2} these mappings can be
uniquely extended to
*-isomorphisms $\P:\eusm{T}\to \A(S')$, $\P_n:\eusm{T}\to \eusm{T}$. Obviously,
$\P_n(S)\to \P(S)$ and by Banach-Steinhaus theorem we have
$\P_n(A)\to \P(A)$ for every $A\in \eusm{T}$. Thus $\P$ is a
nonsurjective *-endomorphism of the Toeplitz algebra which is the
pointwise limit of a sequence of *-automorphisms.

Now, one of the basic results of BDF theory
says that the extensions of $\C(\H)$ by $C(\Bbb T)$ are
completely determined (up to "extension"-equivalence) by the Fredholm
index of
those elements which correspond to the identical function on $\Bbb T$
(see, for example, \cite{Dav, Sections IX.2, IX.3} or \cite{Dou, Section
7}). The extension with Fredholm index $2\leq n\in \N$ can be realized
as the $n$-fold sum of the Toeplitz extension with itself. For instance,
if $n=2$, then this is the algebra of all matrices
$$
\left[
\matrix
K_{11}+A & K_{12}\\
K_{21} & K_{22}+ A
\endmatrix
\right ],
$$
where $A\in \eusm{T}$ and $K_{ij}\in \C(\H)$.
If $\P_n$ is the same as in the case of the Toeplitz algebra above, then
the sequence
$$
\left[
\matrix
K_{11}+A & K_{12}\\
K_{21} & K_{22}+ A
\endmatrix
\right ]
\longmapsto
\left[
\matrix
\P_n(K_{11}+A) & \P_n(K_{12})\\
\P_n(K_{21}) & \P_n(K_{22}+A)
\endmatrix
\right ]
$$
of *-automorphisms converges
pointwise to a nonsurjective *-endomorphism of this algebra. Finally,
since the extensions corresponding to the negative Fredholm index $-n$
are $C^*$-isomorphic (but not "extension"-isomorphic) to the ones with
index $n$, the proof is complete.
\enddemo

\head
3. Remarks and Open Problems
\endhead

We go back to the algebras $\C(\H)+\B(\{ \H_i\})$.
Theorem 1.1 can be considered as an affirmative answer to our
reflexivity problem in the case when there are finitely
many subspaces $\{ \H_i\}$. We feel that it is a natural question to
investigate the same problem for an arbitrary infinite sequence $\{
\H_n\}_{n\in \N}$ of pairwise orthogonal subspaces generating $\H$.
It may be not surprising that this question seems to be
considerably more difficult than the one concerning the finite case.

In \cite{BaMo, Theorem 5} we showed an example of a
nonsurjective *-endomorphism of $\ell_\infty$ which is an
approximately local *-automorphism thus proving the topological
nonreflexivity of the automorphism
and isometry groups of $\ell_\infty$.
This example was based on the existence of a character $\chi$ of
$\ell_\infty$
annihilating $c_0$ which follows from the well-known fact that
$\ell_\infty$ is isomorphic to $C(\beta \N)$, the function algebra on
the Stone-\v Cech compactification $\beta \N$ of $\N$. Suppose now for a
moment that our subspaces $\H_n$ are all one-dimensional. Keeping
the example mentioned above in mind, it is apparent to think of the
map
$$
K+D\longmapsto S(K+D)S^*+\chi (D)P,
$$
where $S$ is a unilateral shift, $P=I-SS^*$, $K\in \C(\H)$, $D\in \B(\{
\H_n\})$ and $\chi(D)$ is the value of $\chi$ on the sequence of
eigenvalues of $D$ corresponding to the eigensubspaces $\H_n$.
Similarly to the argument followed in the proof of \cite{BaMo, Theorem
5} it is easy to see that the above map is an approximately local
*-automorphism of
$\C(\H)+\B(\{ \H_n\})$ which is not surjective. Hence, in this case we
obtain that the automorphism and isometry groups are topologically
nonreflexive.
We note that the previous approach can be
generalized to the case $\sup \dim \H_n <\infty$ quite easily.
In fact, this follows in the same way after referring to the property
of the Stone-\v Cech compactification that not only the complex valued
bounded functions on $\N$ can be uniquely extended to a continuous
function on $\beta \N$, but the same holds true for functions on $\N$
which take their values in a compact Hausdorff space.

Though we know the answer to our problem in the previous particular
case, the question is open in its full generality and reads as
follows.

\proclaim{Problem 3.1}
Let the number of the subspaces $\{ \H_i\}$ be infinite. Are
the automorphism and isometry groups of $\C(\H)+\B(\{ \H_i\})$
topologically nonreflexive? If it is not true in general, then determine
those cases when we have the topological reflexivity.
\endproclaim

The next problem which comes naturally is concerned with
the algebraic reflexivity. The answer to this question is missing even
in
the case $\dim \H_n=1$ $(n\in \N)$. Let us point out to the difficulties
which one has to face.
Since under the isomorphism between $\ell_\infty$ and $C(\beta \N)$,
$c_0$ corresponds to those elements
of $C(\beta \N)$ which vanish on $\N^*=\beta \N \setminus \N$, therefore
in our case $\C(\H)+\B(\{ \H_n\})$ can be considered as an
extension of $\C(\H)$ by $C(\N^*)$. When trying to answer the problem,
one might tend to apply the same argument that we have followed above
concerning extensions.
Unfortunately, for the commutative $C^*$-algebra $C(\N^*)$ we do not
have a result like Proposition 2.8 which played basic role
in the proof of Theorem 1.2.
We note that to the validity of the proof of \cite{MoZa, Theorem 2.2}
and hence to the validity of that of Proposition 2.8,
the only property that we have to require from the underlying compact
Hausdorff space $X$ of $C(X)$ is that to
every point $x\in X$, let there be a continuous complex valued function
$f$ on $X$ for which $f(x)\neq f(y)$ whenever $x\neq y\in X$. But this
property fails for $\N^*$. In fact, in any topological space the
points where a given continuous function takes a given value form a
$G_\delta$-set. But in $\N^*$ every nonempty $G_\delta$-set has nonempty
interior
\cite{Wal, p. 78} which, together with the fact that no point of
$\N^*$ is isolated \cite{Wal, p. 74} give us that every continuous
function takes any value from its range infinitely many times.
Nevertheless, one can attempt to answer the following problem
which seems to be a completely topological question.

\proclaim{Problem 3.2} Are the automorphism and isometry groups of
$C(\N^*)$ algebraically reflexive?
\endproclaim

At the first glance, having the fact in mind that the topological
structure of $\beta \N$ is at least so complicated as that of $\N^*$, it
might be interesting to observe that the
same question for $C(\beta \N)$ is easy to answer (which is "yes", of
course) \cite{MoZa, Theorem 2.1}. This surprise disappears at once if we
take into consideration the main
difference between $\beta \N$ and $\N^*$ which plays significant role
here. Namely, the points of $\N$ are isolated in $\beta \N$, and
thus every homeomorphism of $\beta \N$ leaves $\N$ invariant. One should
compare this with the fact that in $\N^*$ the orbit of every
point has at least $\frak c$ points \cite{Wal, pp. 76-77}.

Since even an affirmative answer to Problem 3.2 would not give a
complete solution of the question of algebraic reflexivity of our groups
in the case of the algebra $\C(\H)+\B(\{ \H_i\})$,
we can formulate our last open problem as follows.

\proclaim{Problem 3.3}
Let the number of the subspaces $\{ \H_i\}$ be infinite. Are the
automorphism and isometry groups of the $C^*$-algebra $\C(\H)+\B(\{ \H_i\})$
algebraically reflexive? If it is not true in general, then
characterize those cases when we have the algebraic reflexivity.
\endproclaim

\enddocument
\end